\documentclass[11pt]{article}

\input{xypic.tex} \xyoption{all}
\usepackage{amsthm,amssymb,amsmath,amscd,latexsym}
\usepackage{graphics}

\newcommand{\real}{\mathbb{R}} \newcommand{\cpx}{\mathbb{C}} 
\newcommand{\zed}{\mathbb{Z}}  \newcommand{\bT}{\mathbb{T}} 
\newcommand{\id}{\mathbf{1}} \newcommand{\ch}{{\rm ch}}
\newcommand{\hE}{\hat{E}} \newcommand{\hA}{\hat{A}}
\newcommand{\cA}{{\cal A}} \newcommand{\cG}{{\cal G}}
\newcommand{\del}{\overline{\partial}}
\newcommand{\p}{\mathbb{P}} \newcommand{\bB}{\mathbb{B}}

\newtheorem{theorem}{Theorem}
\newtheorem{proposition}[theorem]{Proposition}
\newtheorem{lemma}[theorem]{Lemma}

\newtheorem*{remark}{Remark}

\begin{document}

\title{A survey on Nahm transform}
\author{ Marcos Jardim \\ University of Massachusetts at Amherst \\
Department of Mathematics and Statistics \\
Amherst, MA 01003-9305 USA }

\maketitle

\begin{abstract}
We review the construction known as the Nahm transform in a generalized context, which includes
all the examples of this construction already described in the literature. The Nahm transform for
translation invariant instantons on $\real^4$ is presented in an uniform manner. We also analyze two
new examples, the first of which being the first example involving a four-manifold that is not hyperk"ahler.
\end{abstract}


\baselineskip18pt


\section{Introduction}

Since the appearance of the Yang-Mills equation on the mathematical
scene in the late 70's, its anti-self-dual (ASD) solutions have
been intensively studied. The first major result in the field was the
ADHM construction of instantons on $\real^4$ \cite{ADHM}. Soon
after that, W. Nahm adapted the ADHM construction to obtain the
{\em time-invariant} ASD solutions of the Yang-Mills equations, the
so-called monopoles \cite{N1,N}. Nahm found a correspondence between
solutions of the anti-self-duality equations which are invariant under
translations in one direction and solutions of the anti-self-duality equations
which are invariant under translations in three directions. His physical arguments
were formalized in a beautiful paper by N. Hitchin \cite{H3}.

It was later realized that these constructions are two examples of a much more
general framework. This was first pointed out by Corrigan \& Goddard in \cite{CG},
and further elaborated in papers by Braam \& van Baal \cite{BVB} (who coined the
term ``Nahm transform'') and by Nakajima \cite{Na}.

The Nahm transform was initially conceived as a correspondence between solutions of
the anti-self-duality equations which are invariant under dual subgroups of translations
of $\real^4$, and many such correspondences have been described in the literature
(see Section \ref{ex} below). The first goal of this paper, pursued in Section \ref{nahm}, 
is to show that the transform can be set up in a much larger class of four-manifolds, namely
spin manifolds of non-negative scalar curvature. It can be characterized as nonlinear version
of the Fourier transform, which takes vector bundles provided with anti-self-dual connections
over a 4-dimensional manifold into vector bundles with connections over a {\em dual} manifold.
If further geometric structures are available one can easily show that the transformed connection
satisfies certain natural differential conditions. In particular, if the original manifold admits a hyperk\"ahler
metric, then the transformed connection is a quaternionic instanton.

We then list all instances of the Nahm transform described in the literature, adding two new examples.
The first one, concerning instantons on the four-sphere, is of particular interest, for it involves a  four
dimensional manifold which does not admit a complex structure. Despite this, we show that an interesting
Nahm transform mapping intantons into instantons can be defined. 

This paper is written with a wider audience of in mind, so arguments familiar to experts are presented in detail.
We focus on the mathematical aspects and precise mathematical statements surrounding
the Nahm transform. There is an extensive physical literature relating Nahm transform and fundamental problems
in Physics, like quark confinement in QCD and string dualities. For the reader interested in these issues, we recommend
for instance \cite{GGPB} (among other papers by Pierre van Baal) for the relevance of Nahm transform in QCD
on the lattice and \cite{CK,Di,KS,W} for the relations between Nahm transform and string theory. Another interesting
related topic that is out of this survey is the role of Nahm transform in noncomutative gauge theories \cite{ANS,NeS}.


\section{The Nahm transform} \label{nahm}

Let $(M,g_M)$ be a smooth oriented Riemannian spin 4-manifold with non-negative
scalar curvature ($R_M \geq 0$). For simplicity, we assume that $M$
is compact. We denote by $S^\pm$ the spinor bundles of positive and negative
chirality.

Consider a complex vector bundle $E$ over $M$, and let $A$ be an
anti-self-dual connection on $E$; more precisely, its curvature $F_A$
satifies the following condition:
\begin{equation} \label{asd}
\ast F_A = - F_A
\end{equation}  
where $\ast$ denotes the Hodge star operator.
We also assume that $A$ is {\em 1-irreducible}:
$$ \nabla_A s = 0 ~~ \Rightarrow ~~ s = 0 $$
i.e. there are non covariantly constant sections.

Now let $T$ be a smooth manifold parametrizing a family of (gauge equivalence
classes of) irreducible, anti-self-dual connections on a fixed complex vector
bundle $F\to M$. In other words, each $t\in T$ corresponds to an anti-self-dual
connection $B_t$ on the bundle $F$. Typically, we can think of $T$ as a
(submanifold of a) moduli space of irreducible anti-self-dual connections
on $F\to M$. Note also that the Riemannian metric on $M$ induces a natural
metric on $g_T$ on $T$.

The {\em Nahm transform from $M$ to $T$} is a mechanism that transforms
vector bundles with anti-self-dual connections on $M$ into vector bundles with connections
on $T$. If $T$ parametrizes a family of flat connections over $M$, we will say that the
transform is {\em flat}; otherwise, we will say that the Nahm transform is {\em non-flat}

Let us now describe the transform in detail. 
On the tensor bundle $E\otimes F$, we have a twisted family of 
anti-self-dual connections $A_t = A\otimes\id_F+\id_E\otimes B_t$.
So we can consider the family of coupled Dirac operators:
$$ D_{A_t} : L^2_p(E\otimes F\otimes S^+) \longrightarrow
L^2_{p-1}(E\otimes F\otimes S^-) $$
let $D_{A_t}^*$ denote the dual Dirac operator. The Dirac laplacian
$D_{A_t}^*D_{A_t}$ is related to the trace laplacian $\nabla_{A_t}^*\nabla_{A_t}$
via the Weitzenb\"ock formula:
\begin{equation} \label{w1}
D_{A_t}^*D_{A_t} = \nabla_{A_t}^*\nabla_{A_t} - F_{A_t}^+ + \frac{1}{4}R_M
\end{equation}
Applying (\ref{w1}) to a section $s\in L^2_p(E\otimes F\otimes S^+)$, and
integrating by parts, we obtain:
\begin{equation} \label{w2}
|| D_{A_t} s ||^2 = || \nabla_{A_t} s ||^2 + 
\frac{1}{4} \int_M R_M \langle s , s \rangle \geq 0
\end{equation}
with equality if and only if $s=0$, since $F_{A_t}^+=0$ and $R_M\geq 0$.
Therefore, we conclude that $\ker D_{A_t} = \{0\} $ for all $t\in T$.

This means that $\hE = -{\rm Index}\{D_{A_t}\}$ is a well-defined
(hermitian) vector bundle over $T$; the fibre $\hE_t$ is given by
${\rm coker}~D_{A_t}$.

Furthermore, letting $\hat{H}$ denote the trivial Hilbert
bundle over $T$ with fibres given by $L^2_{p-1}(E\otimes F\otimes S^-)$, 
one can also define a connection $\hA$ via the {\em projection formula}:
\begin{equation} \label{proj}
\nabla_{\hA} = P \underline{d} \iota
\end{equation}
where $\iota: \hE\to\hat{H}$ denotes the natural inclusion, $\underline{d}$
denotes the trivial covariant derivative on $\hat{H}$ and $P:\hat{H}\to\hE$
denotes the orthogonal projection induced by the $L^2$ inner product;
at each $t\in T$, this projection can be expressed in the following way:
\begin{equation} \label{projector}
P(t) = \id_{\hat H} - D_{A_t} G_{A_t} D_{A_t}^*
\end{equation}
where $G_{A_t}=\left(D_{A_t}^*D_{A_t}\right)^{-1}$ is the Green's operator
for the Dirac laplacian.

Notice that if  $t,t'\in T$ are such that the corresponding connections
$B_t$ and $B_{t'}$ are gauge equivalent, then clearly $A_t$ and $A_{t'}$
are also gauge equivalent. Hence there is a natural isomorphism
$ \ker D_{A_t}^* \stackrel{\sim}{\rightarrow} \ker D_{A_{t'}}^*$, and 
the index bundle $\hE$ descends to a bundle on the quotient $T/\cG$, where
$\cG$ denotes the group of gauge transformations of $F$. For this reason, we
assume from now on that $T$ parametrizes a family of  gauge equivalence
classes of irreducible anti-self-dual connections on $F$.

The pair $(\hE,\hA)$ is called the {\em Nahm transform} of $(E,A)$.
Clearly, the transformed connection $\hA$ is unitary.

\begin{remark} \rm
The key necessary and sufficient condition for the transform to work is
the vanishing of the kernel of the Dirac operators $D_{A_t}$ for all $t\in T$.
This means that the non-negativity condition on the scalar curvature $R_M$ can be weakened.
Indeed, consider the following bilinear hermitian pairing on $L^2(E\otimes F\otimes S^+)$:
$$ \{s_1,s_2\} := \int_M R_M \langle s_1 , s_2 \rangle~,~~~~
s_1,s_2\in L^2(E\otimes F\otimes S^+) $$
Using theWeitzenb\"ock formula (\ref{w1}), it is easy to see that
$\ker D_{A_t} = 0$ if and only if $\{s,s\}\geq-4||\nabla_{A_t}s||^2$
for all $s\in L^2(E\otimes F\otimes S^+)$ and all $t\in T$, with
equality if and only if $s=0$. 
\end{remark}

\begin{lemma}
If $A$ and $A'$ are two gauge equivalent connections on a vector bundle
$E\to X$, then $\hA$ and $\hA'$ are gauge equivalent connections
on the transformed bundle $\hE\to Y$.
\end{lemma}

In other words, the Nahm transform yields a well-defined map from the moduli space
of gauge equivalence classes of anti-self-dual connections on $E\to M$ into the space
of gauge equivalence classes of connections on $\hE\to T$.

\begin{proof}
Since $A$ and $A'$ are gauge equivalent, there is a bundle automorphism
$h : E\to E$ such that $\nabla_A' = h^{-1} \nabla_A h$. Take
$g=h\otimes{\mathbf 1}_F \in {\rm Aut}(E\otimes F)$, so that
$\nabla_{A_t'} = g^{-1} \nabla_{A_t} g$, hence
$D_{A_t'}^* = g^{-1} D_{A_t}^* g$, for all $t\in T$. Thus if
$\{\Psi_i\}$ is a basis for $\ker D_{A_t}^*$, then 
$\{\Psi_i'=g^{-1}\Psi\}$ is a basis for $\ker D_{A_t'}^*$.
So $g$ can also be regarded as an automorphism of the transformed
bundle $\hE$. It is then easy to see that:
$$ \nabla_{\hA'} = P' \underline{d} \iota' = 
\left( g^{-1} P' g \right) \underline{d}\left( g^{-1} \iota' g \right) = 
g^{-1} \nabla_{\hA} g $$
since $\underline{d}g^{-1}=0$, for $g=h\otimes{\mathbf 1}_E$ does not
depend on $t$.
\end{proof}

The Nahm transformed connection $\hat A$ was defined above in a rather
coordinate-free manner. For many calculations, it is important to have a
more explicit description. First note that the rank of the transformed
bundle $\hE$ is just the index of the Dirac operator $D_{A_t}^*$ for some
$t\in T$, so it is given by:
\begin{equation} \label{rank}
\hat{r} = {\rm rank}~\hE =
- \int_M \ch(E) \cdot\ch(F) \cdot \left( 1 - \frac{1}{24} p_1(M) \right)
\end{equation}
where $p_1(M)$ denotes the first Pontryagin class of $M$. Recall that
since $M$ is a spin 4-manifold, then
$$ p=1/24 p_1(M)[M] = \frac{1}{192\pi^2}\int_M {\rm Tr}(R_M\wedge R_M) $$
is an even integer (so-called $\hat{\mathfrak{a}}$-genus of $M$).

Now let $\{\Psi_i=\Psi_i(x;t)\}_{i=1}^{\hat r}$ be linearly independent
solutions of the Dirac equation $D_{A_\xi}^*\Psi_i=0$. We can assume that
$\langle \Psi_i, \Psi_j \rangle = \delta_{ij}$, where $\langle\cdot,\cdot\rangle$
denotes the $L^2$ inner product on $\hat H$. Clearly,
$\{\Psi_i\}_{i=1}^{\hat r}$ forms a local orthonormal frame for $\hE$.
In this choice of trivialization, the components of the connection matrix
$\hA$ can be written in the following way:
\begin{equation} \label{matrix}
\hA_{ij} = \langle \Psi_i, \nabla_{\hA} \Psi_j \rangle
= \langle \Psi_i, \underline{d} \Psi_j \rangle = 
\int_M \Psi_i(x;t)^\dagger \bullet  \frac{d}{dt} \Psi_j(x;t) ~~ d^4x 
\end{equation}
where $\bullet$ denotes Clifford multiplication.

In this trivialization, the curvature can be expressed as follows:
\begin{eqnarray*}
(F_{\hA})_{ij} & = & \langle \Psi_i, \nabla_{\hA} \nabla_{\hA} \Psi_j \rangle
= \langle \Psi_i, \underline{d} P \underline{d} \Psi_j \rangle = \\
& = & \langle \Psi_i, \underline{d} D_{A_t} G_{A_t} D_{A_t}^* \underline{d} \Psi_j \rangle = 
-\langle D_{A_t}^* \underline{d} \Psi_i, G_{A_t} D_{A_t}^* \underline{d} \Psi_j \rangle 
\end{eqnarray*}
We define $\Delta=[D_{A_t}^*,\underline{d}]$; this is an algebraic operator
acting as:
$$ \Delta: L^2(M\times T,\pi_1^*(E\otimes F \otimes S^-))
\longrightarrow L^2(M\times T,\pi_1^*(E\otimes F \otimes S^+)\otimes\pi_2^*\Omega^1_T) $$
where $\pi_1$ and $\pi_2$ are the projections of $M\times T$ onto the
first and second factors, respectively. More precisely, this operator
can be expressed in terms of Clifford multiplication; in local coordinates:
$$ \Delta = \sum_{k=1}^{\dim T} \delta_k(x;t) dt_k $$
where $\delta_k(x;t)$ is a local section of $\pi_1^*(E\otimes F \otimes S^-)$.
With this in mind, we conclude that:
$$ \Delta(\Psi) = \sum_{k=1}^{\dim T} \delta_k(x;t)\bullet\Psi dt_k = \Delta\bullet\Psi $$
where $\bullet$ denotes Clifford multiplication. Clearly, if $\Psi\in\ker D_{A_t}^*$,
then $D_{A_t}^*\underline{d}\Psi=\Delta\bullet\Psi$,  Therefore, we have:
\begin{equation} \label{curv}
(F_{\hA})_{ij} = 
-\langle \Delta \bullet \Psi_i, G_{A_t} \left(\Delta \bullet \Psi_j\right) \rangle
\end{equation}

It is important to note that the transformed connection $\hA$ is smooth, but
since the parameter space $T$ might not be compact, $\hA$ might not have
finite $L^2$-norm (i.e. finite Yang-Mills action).


\subsection{The topology of the transformed bundle.}
Let us now study the topological invariants of the transformed bundle. Recall
that one can define a {\em universal bundle with connection} over the product
$M\times T$ in the following way \cite{AS}. Let $\cA$ denote the set of all
connections on $F$, and let $\cG$ denote the group of gauge transformations
(i.e. bundle automorphims). Moreover, let $G$ denote the structure group of $F$,
so that $F$ can be associated with a principal $G_E$ bundle $P$ over $M$ by
means of some representation $\rho:G\to\cpx^n$, where $n={\rm rank}~F$. $\cG$ acts
on $F\times\cA$ by $g(p,A)=(g(p),g(A))$; This action has no fixed points, and
it yields a principal $\cG$-bundle $E\times\cA \rightarrow {\cal Q}$,
where ${\cal Q}=E\times\cA/\cG$.

The structure group $G$ also acts on $E\times\cA$, and since this action
commutes with the one by $\cG$, $G$ acts on $\cal Q$. Moreover, the $G$-action
on ${\cal Q}^{\rm ir}=E\times\cA^{\rm ir}/\cG$ has no fixed points, where
$\cA^{\rm ir}$ denotes the set of irreducible connections on $F$. We end up
with a principal $G$ bundle ${\cal Q}^{\rm ir}\rightarrow M\times(\cA^{\rm ir}/\cG)$,
and we denote by $\tilde{\p}$ the associated vector bundle 
${\cal Q}^{\rm ir}\times_\rho\cpx^n$. Since $T$ is a submanifold of 
$\cA^{\rm ir}/\cG$, we define the {\em Poincar\'e bundle} $\p\to M\times T$
as the restriction of $\tilde{\p}$.

The  principal $G$ bundle ${\cal Q}^{\rm ir}$ also has a natural connection
$\tilde{\omega}$, constructed as follows. The space $E\times\cA^{\rm ir}$ has
a Riemannian metric which is equivariant under $G\times\cG$, so that it descends
to a $G$-equivariant metric on ${\cal Q}^{\rm ir}$. The orthogonal complements
to the orbits of $G$ yields the connection $\tilde{\omega}$. Passing to the
associated vector bundle $\tilde{\p}$ and restricting it to $M\times T$ gives
a connection $\omega$ on the {\em Poincar\'e bundle} $\p$.
The pair $(\p,\omega)$ is universal in the sense that 
$(\p,\omega)|_{M\times\{t\}}\simeq(F,B_t)$ \cite{AS}.

The Atiyah-Singer index theorem for families allows us to compute the
Chern characther of the transformed bundle via the formula:
\begin{equation} \label{ch}
\ch{\hE} = - \int_M \ch(E) \cdot \ch(\p) \cdot \left( 1 - \frac{1}{24} p_1(M) \right)
\end{equation}
where the minus sign is needed because $\hE$ is the bundle of cokernels.
The curvature $\Omega$ of the Poincar\'e connection $\omega$ can be easily
computed, see \cite{AS}. In examples, that can then be used to compute
the Chern character of $\p$.


\subsection{Differential properties of transformed connection.}
Since the expression (\ref{curv}) for the curvature of the transformed
connection does not depend explicitly on the curvature of the original
connection $A$, it is in general very hard to characterize any particular
properties of $F_{\hA}$. 

For instance, when the parameter space $T$ is 4-dimensional, one would
like to know whether $F_{\hA}$ is anti-self-dual. This seems to be a
very hard question in general; we now offer a few positive results.

First, note that the algebraic operator $\Delta=[D_{A_t},\underline{d}]$
can also be thought as a section of the bundle $\pi_1^*{\cal L}\otimes\pi_2^*\Omega^1_T$,
where ${\cal L}={\rm End}(E\otimes F\otimes S^-)$. 

\begin{proposition} \label{commut}
If $[G_{A_t},\Delta]=0$, then $F_{\hA}$ is proportional to $\Delta\wedge\Delta$
as a 2-form over the parameter space $T$. In particular, if $T$ is 4-dimensional,
$F_{\hA}$ is anti-self-dual if and only if $\Delta\wedge\Delta$ is a section of
$\pi_1^*{\cal L}\otimes\pi_2^*\Omega^{2,-}_T$.
\end{proposition}
\begin{proof}
If $G_{A_t}\Delta=\Delta G_{A_t}$, it follows from (\ref{curv}) that:
$$ (F_{\hA})_{ij} = 
-\langle \Delta \bullet \Psi_i, \Delta \bullet (G_{A_t}\Psi_j) \rangle =
- \langle \Delta \bullet \Delta \bullet \Psi_i, G_{A_t}\Psi_j \rangle $$
It is then easy to see from the last expression that each component
$(F_{\hA})_{ij}$ is proportional to $\Delta\wedge\Delta$ as a 2-form
over $T$.
\end{proof}

When $M$ is a K\"ahler or hyperk\"ahler manifold, complex analytic methods
can also be useful. We turn to two well-known results concerning these cases.

\begin{proposition} \label{K}
If $M$ and $T$ are K\"ahler manifolds, then the transformed bundle $\hE$
has a natural complex structure, which is compatible with $\hA$. In
particular, the curvature of the transformed connection is of type $(1,1)$.
\end{proposition}

It is imporant to recall that  if $M$ is a K\"ahler manifold, then all 
connected components of the moduli space of anti-self-dual
connections on $M$ are also K\"ahler. We include an outline of the
proof of this well-known result for the sake of completeness, and for
the convenience of the reader.
 
\begin{proof}
The anti-self-dual connection $A_t$ induces a holomorphic structure 
on the tensor bundle $E\otimes F$, and the Dirac operators can be written
in terms of the Dolbeault operators in the following manner:
$$ D_{A_t} = 2 \left( \del_{A_t} - \del_{A_t}^* \right) ~~~{\rm and}~~~
D_{A_t}^* = 2 \left( \del_{A_t}^* - \del_{A_t} \right) $$
Therefore Hodge theory gives identifications for each $t\in T$:
$$ \ker D_{A_t} = \ker\del_{A_t} \oplus \ker \del_{A_t}^* =
H^0(M,E\otimes F) \oplus H^2(M,E\otimes F) $$
$$ \ker D_{A_t}^* = \ker \del_{A_t}^*\cap \ker\del_{A_t} = H^1(M,E\otimes F) $$  
This means that $\hE$ can be identified (as a smooth vector bundle) with the
cohomology of the family Dolbeault complex:
$$ E\otimes F \stackrel{\del_{A_t}}{\longrightarrow} E\otimes F\otimes\Omega^{0,1}_M
\stackrel{\del_{A_t}}{\longrightarrow} E\otimes F\otimes\Omega^{0,2}_M $$
General theory \cite[p. 79-80]{DK} then implies that $\hE$ also has a holomorphic
structure, with which the connection $\hA$ defined via the projection formula
(\ref{proj}) is compatible.
\end{proof}

Recall that a Riemannian 4-manifold $M$ is said to be {\em hyperk\"ahler} if its holonomy
group is contained in $Sp(1)$. This implies that $M$ carries three almost complex structures
$(I,J,K)$ which are parallel with respect to the Levi-Civita connection and satisfy quaternionic
relations $IJ=-JI=K$.

A {\em quaternionic instanton} is a connection $A$ on a complex vector bundle $V$ over a
hyperk\"ahler manifold $T$ whose curvature $F_A$ is of type (1,1) with respect to all
complex structures \cite{BBH2}. In particular, if $T$ is 4-dimensional then a quaternionic
instanton is just an anti-self-dual connection.

\begin{proposition} \label{hK}
If $M$ and $T$ are hyperk\"ahler manifolds, then the transformed connection is
a quaternionic instanton. In particular, if $T$ is 4-dimensional then $\hA$ is
anti-self-dual.
\end{proposition} 

As in Proposition \ref{K}, the hypothesis here is slightly redundant, for if $M$ is hyperk\"ahler, then
all connected components of the moduli space of anti-self-dual connections on $M$ are also
hyperk\"ahler.

\begin{proof}
Each choice of a K\"ahler structure on $M$ induces a choice of a K\"ahler
structure on $T$; by Proposition \ref{K}, $F_{\hA}$ is of type $(1,1)$ with
respect to this structure. Thus $F_{\hA}$ is of type $(1,1)$ with respect to
all K\"ahler structures on $T$, which means that $\hA$ is a quaternionic instanton.
\end{proof}

Since the only compact 4-dimensional hyperk\"ahler manifolds are
the 4-torus and the K3-surface, this last result seems to have
a rather limited applicability. However, as we will argue in Section
\ref{ex} below, Proposition \ref{hK} can also be used to define a
Nahm transform for instantons over hyperk\"ahler ALE spaces and over
the 4-sphere.

It is also important to mention that a higher dimensional generalization of the Nahm transform for
quaternionic instantons over hyperk\"ahler manifolds has been described by Bartocci, Bruzzo,
and Hern\'andez Ruip\'erez \cite{BBH2}.

\begin{remark} \rm
Finally, we would like to notice that the construction here presented is essentially {\em topological}, 
in the sense that its main ingredient is simply index theory. All the geometric structures used in Section
\ref{nahm} (spin structure, positivity of scalar curvature, hyperk\"ahler metric, etc.) were needed either
because a particular differential operator was used (i.e. the Dirac operator), or
because we selected those objects (i.e. anti-self-dual connection over hyperk\"ahler
manifolds) that yielded very particular transforms (anti-self-dual connections).

One can conceive, for instance, a similar construction either based on a different pseudodifferential
elliptic operator, other than the Dirac operator, or allowing for classes in $K(T)$, rather than actual
vector bundles over the parameter space. The author thus believes that a much more
general construction in a ``K-theory with connections'', akin to the Fourier-Mukai transform
in the derived category of coherent sheaves over algebraic varieties, underlies the construction
here presented. We hope to address this issue in a future paper.
\end{remark}


\section{Examples} \label{ex}

As we mentioned in the Introduction, several examples of the Nahm transform have been
described in the literature, and we now take some time to revise them. 

\subsection{Invariant instantons on $\mathbf{\real^4}$ \& dimensional reduction}

First, we consider the case of translation invariant instantons on $\real^4$, for which the
Nahm transform was first developed. Let $\Lambda$ be a subgroup of translations $\real^4$; the dual group
$$ \Lambda^*=\{\alpha\in(\real^4)^* ~~ | ~~ \alpha(\lambda)\in\zed ~ \forall\lambda\in\Lambda\} $$
can be regarded as a subgroup of translations $(\real^4)^*$. With this in mind, we set $M=\real^4/\Lambda$,
and $T=(\real^4)^*/\Lambda^*$. 

A point $\xi\in T$ can be cannonically identified with  the flat connection $i\cdot\xi$, with $\xi$ being regarded as a (constant)
1-form on $M$, on a topologically trivial line bundle over $M$. Thus all of the Nahm transforms included in this example are
{\em flat}. Conversely, it is easy to see that a point $x\in M$ can also be thought
as  the flat connection $i\cdot x$ on a topologically trivial line bundle over $T$.

At this point it might be useful to briefly remind the reader of the various gauge theoretical equations obtained from the
anti-self-duality equations via dimensional reduction. A connection on a hermitian vector bundle over $\real^4$ of
rank $n$ can be regarded as 1-form
$$ A = \sum_{k=1}^4 A_k(x_1,\cdots,x_4) dx^k ~~,~~ A_k:\real^4\to\mathfrak{u}(n) $$
Assuming that the connection components $A_k$ are invariant under translation in one direction, say $x_4$,
we can think of $\underline{A}=\sum_{k=1}^3 A_k(x_1,x_2,x_3) dx^k$ as a connection on a hermitian vector bundle
over $\real^3$, with the fourth component $\phi=A_4$ being regarded as a bundle endomorphism (the Higgs field). In this way, the
anti-self-duality equations (\ref{asd}) reduce to the so-called Bogomolny (or monopole) equation:
\begin{equation} \label{bogomolny}
F_{\underline{A}} = \ast d\phi
\end{equation}
where $\ast$ is the euclidean Hodge star in dimension 3.

Now assume that the connection components $A_k$ are invariant under translation in two directions, say $x_3$ and $x_4$.
Consider $\underline{A}=\sum_{k=1}^2 A_k(x_1,x_2) dx^k$ as a connection on a hermitian vector bundle
over $\real^2$, with the third and fourth components combined in a complex bundle endomorphim: $\Phi=(A_3+i\cdot A_4)(dx_1-i\cdot dx_2)$.
The anti-self-duality equations (\ref{asd}) are then reduced to the so-called Hitchin's equations:
\begin{equation} \label{hitchin} \left\{ \begin{array}{l}
F_{\underline{A}} = [\Phi,\Phi^*] \\ \del_{\underline{A}} \Phi = 0
\end{array} \right. \end{equation}

Finally, assume that the connection components $A_k$ are invariant under translation in three directions, say $x_2,x_3$ and $x_4$.
After gauging away the first component $A_1$, the anti-self-duality equations (\ref{asd}) reduce to the so-called Nahm's equations:
\begin{equation} \label{nahm eqn}
\frac{d T_k}{dx_1} + \frac{1}{2}\sum_{j,l} \epsilon_{kjl}[T_j,T_l] = 0 ~,~~ j,k,l=\{2,3,4\}
\end{equation}

Roughly speaking, {\em the Nahm transform yield a 1-1 correspondence between $\Lambda$-invariant instantons on $\real^4$
and $\Lambda^*$-invariant instantons on $(\real^4)^*$}. Except for the case $\Lambda=\zed^4$, both $M$ and $T$ are non-compact.
This case is also the only one that relates instantons to instantons, and does not involve a dimensional reduction on either side of the
correspondence.

There are plenty of examples of the Nahm transform for translation invariant instantons available in
the literature, namely:

\begin{enumerate}
\item The trivial case $\Lambda=\{0\}$ is closely related to the
celebrated ADHM construction of instantons, as described by Donaldson
\& Kronheimer \cite{DK}; in this case, $\Lambda^*=(\real^4)^*$ and
an instanton on $\real^4$ corresponds to some algebraic datum (ADHM datum).

\item $\Lambda=\real$ gives rise to monopoles, extensively studied by
Hitchin \cite{H3}, Donaldson \cite{D}, Hurtubise \& Murray \cite{HM} 
and Nakajima \cite{Na}, among several others; here, $\Lambda^*=\real^3$, 
and the transformed object is, for SU(2) monoples, an analytic solution of
Nahm's equations defined over the open interval $(-1,1)$ and with simple
poles at the end-points.

\item If $\Lambda=\zed^4$, this is the Nahm transform of  Schenk \cite{S}, 
Braam \& van Baal \cite{BVB} and Donaldson \& Kronheimer \cite{DK},
defining a correspondence between instantons over two dual four dimensional tori.

\item $\Lambda=\zed$ correspond to the so-called calorons, studied by Nahm \cite{N},
van Baal \cite{VB} and others (see \cite{Ny} and the references therein); the transformed
object is the solution of Nahm-type equations on a circle.

\item The case $\Lambda=\zed^2$ (doubly-periodic instantons) has been
analyzed in great detail by the author \cite{J1,J2,J3} and Biquard \cite{BiJ}.
here, $\Lambda^*=\zed^2\times\real^2$, and the Nahm transform gives a
correspondence between doubly-periodic instantons and certain {\em tame}
solutions of Hitchin's equations on a 2-torus.

\item $\Lambda=\real\times\zed$ gives rise to the periodic monopoles
considered by Cherkis and Kapustin \cite{CK}; in this case, $\Lambda^*=\zed\times\real$,
and the Nahm dual data is given by certain solutions of Hitchin's equations on a 
cylinder. 
\end{enumerate}

In the following two Sections we will take a closer look at periodic instantons and monopoles.


\subsection{Periodic instantons}
Let us now focus on the case of periodic instantons, that is $\Lambda=\zed^d$ and $M=\bT^d\times\real^{4-d}$,
where $d=1,2,3,4$; in these cases, $\Lambda^*=\zed^d\times\real^{4-d}$ and $T=\hat{\bT}^d$.
Other useful accounts of the Nahm transform for periodic instantons in the physical literature
can be found at \cite{FPTW,GGPB}, for example. 

In all the above examples, the general statement one can prove is that there exists a 1-1 correspondence
between instantons over $M$ and singular solutions of the dimensionally reduced anti-self-duality equations over $T$.

Indeed, the correspondence is established just as explained in the previous Section, with some minor modifications
needed to deal with the non-compactness of $M$. Let $T_F(E,A)$ denote set of all points $\xi\in T=\bT^d$
(regarded as a trivial bundle with flat connection) such that the Dirac operator coupled with the tensor connection
$A_\xi=A\otimes\id+\id\otimes\xi$ is Fredholm. Roughly speaking, $T_F(E,A)$ depends only on the 
asymptotic behaviour of the connection $A$, and not on the topological invariants of the bundle $E$; it consists of
$T$ minus finitely many points.

With this in mind, $T_F(E,A)$ can be regarded as parametrizing a family of elliptic Fredholm operators $D_{A_\xi}$
on the bundles $E \to M$. Given that $M$ is flat as a Riemannian manifold, the Weitzenb\"ock formula
\eqref{w1} can be used to show that $\ker D_{A_\xi}=0$ for all $\xi\in T_F(E,A)$, so that $\hE = -{\rm Index}\{D_{A_t}\}$
is a hermitian vector bundle over $T_F(E,A)$. Now $\hE$ can be lifted to a bundle over (a open subset of) $(\real^4)^*$.
A connection $\hA$ on the lifted bundle is defined via the projection formula \eqref{proj}, and $\hA$ can be seen to be
anti-self-dual via the hyperk\"ahler rotation argument in Proposition \ref{hK}. Now $\hA$ descends to the quotient
$T_F(E,A)$, and thus defines a solution of the dimensionally reduced anti-self-duality equations.
Finally, this procedure is invertible, since $M$ can also be regarded as parametrizing trivial line bundles with flat
connections over $T$. 

This simplified statement is still not proven in full generality; only the compact cases $d=4$ and $d=2$ have been
fully described in the literature. The compact case ($d=4$) is the easiest one, and it is closely related to the celebrated
Fourier-Mukai transform in algebraic geometry; see for instance \cite{BVB,DK}. A precise result in this case is as
follows:
\begin{theorem}
There exists a 1-1 correspondence between the following objects:
\begin{itemize} \label{t4}
\item $SU(n)$ instantons over $M=\bT^4$, of charge $k$;
\item $SU(k)$ instantons over $M=\hat\bT^4$, of charge $n$.
\end{itemize}\end{theorem}

The analysis of the non-compact cases ($d=1,2,3$) involve, as we mentioned above, a careful study of the instanton's
asymptotic behaviour, checking that the coupled Dirac operator is indeed Fredholm and correctly applying the Fredholm
theory. The key issue to understand is how the asymptotic data gets transformed.

Doubly-periodic instantons have been extensively studied by the author in \cite{BiJ,J1,J2,J3}. Here is the full statement
of the correspondence, taking into account the asymptotic behaviour of instantons and the singularities of the transformed Nahm data,
in the simplest case of $SU(2)$ gauge group: 
\begin{theorem} \label{dp}
There exists a 1-1 correspondence between the following objects:
\begin{itemize}
\item An anti-self-dual $SU(2)$ connection $A$ on a rank 2 vector bundle $E\to\bT^2\times\real^2$ such that
$$ \frac{1}{8\pi^2} \int_{\bT^2\times\real^2} |F_A|^2 = k ~,$$
and whose asymptotic expansion, up to gauge transformations, as $r\to\infty$ and for
some $\xi=\lambda_1+i\lambda_2\in\hat{\bT}^2$, $\mu=\mu_1+i\mu_2\in\cpx$, and $\alpha\in [0,1/2)$, is
given either by:
$$ i\left(\begin{array}{cc}a_0 & 0 \\ 0 & - a_0 \end{array}\right) + O(r^{-1-\delta})  ~ , ~ {\rm with}$$
$$  a_0 = \lambda_1 dx + \lambda_2 dy + (\mu_1 \cos\theta - \mu_2 \sin\theta)
\frac{dx}{r} + (\mu_1 \sin\theta + \mu_2\cos\theta) \frac{dy}{r} + \alpha d\theta ~, $$
if $\xi,\mu,\alpha \neq 0$; or, if $\xi,\mu,\alpha = 0$, by: 
$$ i\left(\begin{array}{cc} -1 & 0 \\ 0 & 1 \end{array}\right) \frac{d\theta}{\ln r^2} + \frac{1}{r\ln r^2}
\left(\begin{array}{cc} 0 & -\overline{a_0} \\ a_0 & 0 \end{array}\right)+ O(r^{-1}(\ln r)^{-1-\delta}) ~, $$
$$ {\rm with} ~~ a_0 = -e^{i\theta}(dx+idy) $$
\item An hermitian connection $B$ on a rank $k$ hermitian vector bundle $V\to\hat{\bT}^2\setminus\{\pm\xi\}$ and
a skew-hermitian bundle endomorphism $\Phi$ (the Higgs field) satisfying Hitchin's equations:
$$ \left\{ \begin{array}{l}
F_B = [\Phi,\Phi^*] \\ \del_B\Phi = 0
\end{array} \right. $$
and having at most simple poles at $\pm\xi$. Moreover, the residue of $\Phi$ either has rank one, if $\xi\neq-\xi$,
or has rank two, if $\xi=-\xi$, with $\pm\mu$ being the only nonzero eigenvalues; similarly the monodromy of the
connection $B$ near the punctures is semisimple, with either only one nontrivial eigenvalue $\exp(\mp 2\pi i\alpha)$,
or two if $\xi_0=-\xi_0$.
\end{itemize}\end{theorem}

The main feature of the above statement is the matching of the instanton's asymptotic behaviour with the Nahm transformed data's
singularity behaviour.  

It is certainly possible to generalize this correspondence for higher rank (see \cite{J2}), but that would require a much more
lengthy analysis of both the asymptotic behaviour of $A$ and the singularity data of $(B,\Phi)$. It suffices to say that the
while the instanton number $k$ determines the rank of the Nahm transformed bundle $V$, the rank of the original instanton
$A$ determines the number of poles of the transformed Higgs field $\Phi$ (counted according with the rank of its residues).

One expects similar statements to hold also in the cases $d=1$ (calorons) and $d=3$ (spatially periodic instantons);
although the general features of the Nahm transform in these cases are certainly known \cite{KS,Na,VB}, a complete
statement showing how the instantons asymptotic behaviour gets translated into the singularity behaviour for the Nahm
transformed data is still missing.

Some positive results are available calorons. An $L^2$-index theorem for the Dirac operator coupled to calorons
has been established by Nye and Singer \cite{NS}, while the Nahm transform itself has been studied by Nye in his
thesis \cite{Ny}. Nye has identified the appropriate asymptotic behaviour for calorons, and the corresponding
singularity behaviour for the Nahm data on the dual circle $S^1$. He has also constructed the Nahm transform from
calorons to Nahm data on $S^1$ and from Nahm data on $S^1$ to calorons; however, he has not proved that these are
mutually inverse, something that can probably be done using holomorphic geometry and the cohomological argument of
\cite{CK,DK,J2}.

Morever, it also reasonable to expect that the above results for $d=2,4$ (as well as the expected ones for $d=1,3$)
can be adapted to deal with $\zed_p${\em -equivariant} instantons on $\bT^d\times\real^{n-d}$.


\subsection{Periodic monopoles}
The case of periodic monopoles, that is $\Lambda=\zed^d\times \real$, where $d=0,1,2$. As in the case of instantons,
the Nahm transform yields a correspondence between the following objects:
\begin{itemize}
\item monopoles on $M=\bT^d\times\real^{3-d}$;
\item solutions of the dimensionally reduced anti-self-duality equations over $T=\bT^{d}\times\real$.
\end{itemize}

The non-periodic case ($d=0$) was first described by Hitchin in his classical paper \cite{H3} in the simplest
case of gauge group $SU(2)$, and later generalized by Hurtubise \& Murray \cite{HM} to include all classical groups.

\begin{theorem} \label{m}
There exists a 1-1 correspondence between the following objects:
\begin{itemize}
\item An $SU(2)$ connection $A$ on a rank 2 vector bundle $E\to\real^3$ and a skew-hermitian
bundle endomorphism $\Phi$ (the Higgs field) satisfying the Bogomolny equation (\ref{bogomolny}).
and whose asymptotic expansion as $r\to\infty$ is given by, up to gauge transformations and for
some positive integer $k$ (the monopole number):
$$ \Phi \sim \left( \begin{array}{cc} i & 0 \\ 0 & -i \end{array} \right) \cdot \left( 1 - \frac{k}{2r} \right) + O(r^{-2}) $$
$$ |\nabla_A\Phi| \sim O(r^{-2}) ~~~ {\rm and} ~~~ \frac{\partial |\Phi|}{\partial r} \sim  O(r^{-2}) ~ . $$
\item An hermitian connection $\nabla$ on a rank $k$ hermitian vector bundle $V$ over the open interval $I=(-1,1)$ and
three skew-hermitian bundle endomorphisms $T_a$ ($a=1,2,3$) satisfying Nahm's equations (\ref{nahm eqn}),
and having at most simple poles at $t=\pm1$, but are otherwhise analytic. Moreover, the residues of $(T_1,T_2,T_3)$ define
an irreducible representation of $\mathfrak{su}(2)$ at each pole.
\end{itemize} \end{theorem}

The case of periodic monopoles ($d=1$) is studied by in detail Cherkis \& Kapustin \cite{CK}:

\begin{theorem} \label{per-m}
There exists a 1-1 correspondence between the following objects:
\begin{itemize}
\item An $SU(2)$ connection $A$ on a rank 2 vector bundle $E\to S^1\times\real^2$ and a skew-hermitian
bundle endomorphism $\phi$ (the Higgs field) satisfying the Bogomolny equation (\ref{bogomolny}),
and whose asymptotic expansion as $r=|x|\to\infty$ is given by, up to gauge transformations and for
some positive integer $k$ (the monopole number) and parameters $v,w\in\real$:
$$ A \sim w + \left( \begin{array}{cc} i & 0 \\ 0 & -i \end{array} \right) \cdot \frac{k}{2\pi}\theta + O(r^{-1}) $$
$$ \phi \sim v + \left( \begin{array}{cc} i & 0 \\ 0 & -i \end{array} \right) \cdot \frac{k}{2\pi}\log r + O(1) $$
$$ |\nabla_A\Phi| \sim O(r^{-1}) ~~~ {\rm and} ~~~ \frac{\partial |\Phi|}{\partial r} \sim  O(r^{-2}) ~ . $$
\item An hermitian connection $B$ on a rank $k$ hermitian vector bundle $V\to\hat{S^1}\times\real\simeq\cpx^*$ and
$\Phi$ satisfying Hitchin's equations (\ref{hitchin}), and whose asymptotic expansion as $s\to\infty$ are given by, 
up to gauge transformations:
$$ |F_B| \sim O(|s|^{-3/2}) ~ ,$$
$$ {\rm Tr}(\Phi(s)^\alpha) ~ {\rm is~bounded~for} ~ \alpha=1,2,\cdots,k-1 $$
$$ {\rm and} ~~  \det \Phi(s) \sim e^{-2\pi(v+iw)}\cdot O(e^{\pm 2\pi s}) $$
\end{itemize} \end{theorem}

A careful study of doubly-periodic monopoles (the $d=2$ case) is still lacking. It is interesting to note that the Nahm transform of
doubly-periodic monopoles is {\em self-dual}, in the sense that $M=T=\bT^2\times\real$; in other words, the Nahm transform takes
doubly-periodic monopoles into (singular) doubly-periodic monopoles, probably permutating rank and charge.


\subsection{K3 surfaces}
A very interesting example of a {\em non-flat} Nahm transform was described by Bartocci, Bruzzo and
Hern\'andez-Ruip\'erez in \cite{BBH1,BBH2}. Let $M$ be a reflexive $K3$ surface, which is defined by
the following requirements:
\begin{enumerate}
\item $M$ admits a K\"ahler form $\omega$ whose cohomology class $H$ satisfies $H^2=2$;
\item $M$ admits a holomorphic line bundle $L$ whose Chern class $\ell=c_1(L)$ is such that
$\ell\cdot H=0$ and $\ell^2=-12$;
\item if $D$ is the divisor of a nodal curve on $M$, one has $D\cdot H > 2$.
\end{enumerate}

Now let $T$ be the moduli space of instantons of rank 2 with determinant line bundle $L$ (so that $c_1=\ell$)
and $c_2=-1$ over $M$; it can be shown that $T$ is isomorphic to $M$ as a complex algebraic variety \cite{BBH1}.
Since both $M$ and $T$ are hyperk\"ahler manifolds, Nahm transform takes instantons over $M$ into instantons over
$T$. Furthermore, under appropriate circunstances, the transform is invertible, and one obtains in particular the following
result \cite{BBH1,BBH2}:

\begin{theorem}
There exists a 1-1 correspondence between the following objects ($n\geq2$ and $k\geq1$):
\begin{itemize}
\item $SU(n)$ instantons of charge $k$ over $M$;
\item $U(2n+k)$ instantons of charge $k$ over $T$, with first Chern class given by $k\hat{\ell}$.
\end{itemize}\end{theorem}

Finally, we would like to point out that a similar result also holds for hyperk\"ahler ALE 4-manifolds; a
preliminary version was announced in \cite{BaJ} (see also \cite{GN}).


\subsection{First new example: doubly-periodic instantons.}
Let us now proceed to describe two new examples of {\em non-flat} Nahm transforms.
The second one, described below, is particularly interesting, for it is the only example
in which $M$ is not a hyperk\"ahler 4-manifold.

Our first new example of a non-flat Nahm transform is based on the observation that, once
asymptotic parameters $(\xi,\mu,\alpha)$ are fixed, the moduli space ${\cal M}_{(1,\xi,\mu,\alpha)}$
of charge one $SU(2)$ doubly-periodic instantons (as described in Theorem \ref{dp}) is just
$\bT^2\times\real^2$ with the flat metric \cite{BiJ}.

Thus set $M=\bT^2\times\real^2$ and $T={\cal M}_{(1,\xi,\mu,\alpha)}=\bT^2\times\real^2$;
let $E\to M$ be a hermitian vector bundle of rank $n$, and let $A$ be an anti-self-dual
connection on $E$. Denote the points of $T$ by the pair $(F,B)$ consisting of a rank 2
hermitian vector bundle $F$ and an anti-self-dual connection $B$. If the asymptotic state
of the connection $A$ does not contain $\xi$, then the twisted connection $A_B=A\otimes\id + \id\otimes B$
contains no flat factors at infinity, and the Dirac operators $D_{A_B}^\pm$ are Fredholm
\cite{J2}. This means that the Nahm transformed bundle with connection $(\hat{E},\hat{A})\to T$
are well-defined, according to procedure in Section \ref{nahm}. Using the hyperk\"ahler rotation
method of Proposition \ref{hK}, one sees that $\hat{A}$ is also anti-self-dual.

Clearly, $M$ can also be regarded as a moduli space of instantons on $T$, so there is a Nahm
transform that transforms instantons on $T$ into instantons $M$. It seems reasonable to
conjecture that these transforms are the inverse of one another.


\subsection{Second new example: instantons over the 4-sphere.}

Let us now briefly analyse the Nahm transform for the simplest possible
compact spin 4-manifold with non-negative scalar curvature.
Let $M=S^4$ be the round 4-dimensional sphere, and let $T$
be the moduli space of $SU(2)$ instantons over $S^4$ with charge
one; as a Riemannian manifold, $T$ is a hyperbolic 5-ball $\bB^5$ \cite{FU}.

So let $E\to S^4$ be a complex vector bundle of rank $n\geq2$, provided
with an instanton $A$ of charge $k\geq1$. Nahm transform gives a bundle
$\hE\to B^5$ of rank $2k+r$, by the index formula (\ref{rank}). Since $\bB^5$
is simply-connected, this is the only nontrivial topological invariant of
the transformed bundle. This illustrates the wide range of possibilities for a
Nahm transform beyond the confines of hyperk\"ahler geometry.





\paragraph{Acknowledgement.}
This survey was largely motivated by a talk delivered by Michael Atiyah at the conference
``Unity of Mathematics'' (Harvard University, September 2003). In his talk, Atiyah proposed
that non-linear versions of the Fourier transform should play a crucial role in 21$^{\rm st}$
Mathematics. The author thanks him for that very inspirational talk, and Claudio Bartocci,
Ugo Bruzzo and James Glazebrook for their comments on the preliminary version.

 \end{document}